\documentclass[a4paper,11pt]{amsart}

\usepackage{amsmath,amscd,amsfonts,amsthm,amssymb,eucal}
\usepackage{url}

\usepackage[]{latexsym}

\newtheorem{thm}{Theorem}
\newtheorem{cor}[thm]{Corollary}

\theoremstyle{definition}

\newtheorem*{rem}{Remark}

\newenvironment{pf}{\par\noindent{\bf Proof.}\enspace\ignorespaces}{\qed\par\par}
\def\qed{\relax\ifmmode\hskip2em \Box\else\unskip\nobreak\hskip1em $\Box$\fi}

\newcommand{\bQ}{{\mathbb{Q}}}

\newcommand{\bR}{{\mathbb{R}}}
\newcommand{\bZ}{{\mathbb{Z}}}
\newcommand{\bA}{{\mathbb{A}}}

\newcommand{\Jac}{\mbox{Jac}}

\newcommand{\rank}{\mbox{rank }}

\begin{document}

\title{Hyperelliptic curves covering an elliptic curve twice}
\author{Xavier Xarles}
\address{Departament de Matem\`atiques\\Universitat Aut\`onoma de
Barcelona\\08193 Bellaterra, Barcelona, Catalonia}
\email{xarles@mat.uab.cat}

\thanks{The author was partially supported by the grant
MTM2009-10359.} \maketitle

\begin{flushright}
    \textsc{\textit{To the memory of F. Momose}}
  \end{flushright}

\vspace{1cm}

Let $K$ be a field and let $E$ be an elliptic curve over $K$. A
natural question that have been considered by Mestre in \cite{Me}
is the existence of an hyperelliptic curve $H$ defined over $K$
with two independent maps from $H$ to $E$, hence such that its
Jacobian $\Jac(H)$ is isogenous to $E^2\times A$, for some abelian
variety $A$. Mestre constructs in (op. cit.) such a curve $H$ with
genus 6 (for a field with characteristic 0). In this short note we
show that there exists such a curve, but with genus 5, and for any
field of characteristic $\ne 2,3$. We don't know if the
characteristic $2$ and $3$ cases can be solved using the same
methods. We also show that, if we want such a curve $H$ to be just
geometrically hyperelliptic, so having a degree two map to a
conic, then there is one with genus 3. It is known that there does
not exists such a curve with genus 2 for a general elliptic curve
over a general field.

These results have consequences on the distribution of rank $\ge
2$ twists of elliptic curves over $\bQ$, as showed by Steward and
Top (\cite{ST}, see also \cite{RS}).

We start with the second result. The following theorem was already
been shown (but not in this form) by Mestre in a remark in
\cite{Me} and by Steward and Top in \cite{ST} during the proof of
their Theorem 4. The proof is elementary and we leave it to the
reader.

\begin{thm} Let $E$ be an elliptic curve over $K$ given by a Weierstrass
equation of the form $y^2=x^3-Ax+B$, for some $A\ne 0$ and $B$ in
$K$. Then the curve $H$ obtained form the desingularization of the
projectivization of the curve in the affine space $\bA^3$ given by
the equations
$$ \left. \begin{aligned} y^2=x^3-Ax&+B \\ x^2+xz+z^2&=A \end{aligned}\right\} $$
is geometrically hyperelliptic and it has two independent natural
maps to $E$, given by $f_1(x,y,z)=(x,y)$ and $f_2(x,y,z)=(z,y)$.
\end{thm}

\begin{rem}
It is a short computation to show that the abelian variety
$\Jac(C)$ is isogenous to $E^2\times E'$, where $E'$ is the
elliptic curve given by the Weierstrass equation
$y^2=x^3-27Bx^2+27A^3x$.
\end{rem}

Next theorem is based on a similar idea of the previous theorem,
but using a quartic equation instead of a cubic equation for an
elliptic curve.

\begin{thm} Let $K$ be a field of characteristic $\ne 2$ and $3$,
and take $j\in K$, $j\ne 0$ and $1728$. Consider
$A:=\frac{3^3j}{2^2(j-1728)}$. Then the elliptic curve $E$ over
$K$ given by the Weierstrass equation $y^2=x^3-Ax+A$ has
$j$-invariant $j$ and the hyperelliptic curve $H$ of genus $5$
given by the equation
$$y^2=A(x+1)^4(x^2+1)^4-2^6x^3(x^2+x+1)^3$$
has two independent maps to $E$.
\end{thm}

\begin{pf}
The assertion on the $j$-invariant is easy by using the known
formulae of the $j$-invariant. We are going to find the two maps
from $H$ to $E$ by writing $H$ in a different form.

First of all, we write the genus $1$ curve $E$ as the curve $D$
with equation of the form $y^2=x^4+x^3+B$, where $B=\frac A{2^6}$.
It is standard to show that both curves are isomorphic. Now,
consider the affine (singular) curve given by the equation
$$F(x,z)=\frac{(x^4-z^4)}{(x-z)}+\frac{(x^3-z^3)}{(x-z)}=0$$
in the affine plane $\bA^2$. This curve has genus $0$ and can be
parametrized by
$$ x=-\frac{t^3-1}{t^4-1} \ , \ z=-t\frac{t^3-1}{t^4-1}$$
By using this formulae we get that the curve $C$ obtained form the
desingularization of the projectivization of the curve in the affine
space $\bA^3$ given by the equations
$$ \left. \begin{aligned} y^2=x^4+x^3&+B \\
\frac{(x^4-z^4)}{(x-z)}+&\frac{(x^3-z^3)}{(x-z)}=0 \end{aligned}\right\}$$
is isomorphic to the hyperelliptic curve $H$ given by the equation
$$y^2=A(x+1)^4(x^2+1)^4-2^6x^3(x^2+x+1)^3.$$
Now, the two independent maps from $H$ to $E$ are described easily
as maps from $C$ to $D$ given by $f_1(x,y,z)=(x,y)$ and
$f_2(x,y,z)=(z,y)$. Both maps are clearly well defined. To show
they are independent, it is sufficient to prove that there does
not exist two distinct endomorphisms $m_1$ and $m_2$ of $E$ as
genus 1 curve such that $m_1 \circ f_1=m_2\circ f_2$. First,
recall that $f_1$ and $f_2$ send the only point at infinity to the
$0$ point of $E$, hence $m_1$ and $m_2$ must be endomorphisms as
elliptic curves. Moreover, both $f_1$ and $f_2$ have the same
degree (equal to 3), hence the only possibility is that
$f_1=-f_2$, which is clearly not the case. \end{pf}

\begin{rem}
It's not difficult to show that the curve $H$ is an unramified
degree $2$ covering of the genus $3$ curve $H_1$ given by the
equation
$$y^2=(x^2-4)\left(A(x+2)^2x^4-2^6(x+1)^3\right)$$ and it has an
independent degree $2$ covering to the genus $2$ curve $H_2$ given
by the equation $y^2=A(x+2)^2x^4-2^6(x+1)^3$. The curve $H_2$ has
degree $2$ maps to the curve $E$ defined above, and the curve $E'$
given by the Weierstrass equation $y^2=x^3+Ax^2+2Ax+A$. The curve
$H_1$ has a map to the curve $E$, and its kernel is simple in
general (a fact that can be shown for some value of $A$ and modulo
some prime $p$ by computing its zeta function). Hence the abelian
variety $\Jac(H)$ is isogenous to $E^2\times E'\times F$, where
$F$ is an abelian surface.
\end{rem}

For any (hyper)elliptic curve $E$ defined over a field $K$ given
by an equation of the form $y^2=p(x)$, where $p(x)\in K[X]$, and
for any $d\in K^*$, we will denote by $E_d$ the quadratic twist of
$E$, which is the (hyper)elliptic curve given by the equation
$dy^2=p(x)$.

\begin{cor} Let $K$ be a field with characteristic $\ne 2,3$, and let
$E$ be an elliptic curve over $K$. Then there exists an
hyperelliptic curve $H$ of genus $\le 5$ and defined over $K$ with
two independent maps from $H$ to $E$ defined over $K$.
\end{cor}
\begin{pf} If the curve $E$ has $j$-invariant equal to $0$ or
$1728$, the result is well known (and in fact one can find such a
curve with genus $2$, see for example \cite{ST}).

Now, if the $j$-invariant is $\ne 0$ and $1728$, there exists a
quadratic twist of $E$ isomorphic to the elliptic curve given by
the equation $y^2=x^3-Ax+A$, for $A:=\frac{3^3j}{2^2(j-1728)}$.
Then the same quadratic twist of the curve $H$ given by the
previous theorem gives the answer. \end{pf}

We can apply this result to give lower bounds for the number of
quadratic twists with rank at least 2 for elliptic curves over
$\bQ$, improving slightly the results by Steward and Top in
\cite{ST}. The result is a direct application of their theorem 2,
but using our theorem 2 instead of Mestre's construction.

\begin{cor} Let $E$ be an elliptic curve over $\bQ$. Consider, for
any positive real number $X\in \bR$, the function
$$N_{\ge 2}(X):=\#\{ \mbox{squarefree } d \in \bZ : |d|\le X \mbox{
and } \rank_{\bZ}(E_d(\bQ))\le 2\}.$$ Then $N_{\ge 2}(X)\gg
X^{1/6}/\log^2(X)$.
\end{cor}

\end{document}